\documentclass[12pt]{amsart}
\usepackage{graphicx,amssymb}
\usepackage[all]{xy}
\newtheorem{thm}{Theorem}[section]
\newtheorem{cor}[thm]{Corollary}
\newtheorem{lem}[thm]{Lemma}
\newtheorem{prop}[thm]{Proposition}
\theoremstyle{definition}

\newtheorem{rem}[thm]{Remark}

\numberwithin{equation}{section}

\newcommand{\cD}{\mathcal{D}}
\newcommand{\ra}{\rightarrow}
\def\lra{\longrightarrow}
\def\bs{{\backslash}}
\DeclareMathOperator{\im}{Im}
\def\P{{\mathbb P}}

\begin{document}

\title[ ]{Products of Jacobians as Prym-Tyurin varieties}%

\author{A. Carocca, H. Lange, R. E. Rodr{\'\i}guez and A. M. Rojas}

\address{A. Carocca\\Facultad de Matem\'aticas,
Pontificia Universidad Cat\'olica de Chile, Casilla 306-22,
Santiago, Chile}
\email{acarocca@mat.puc.cl}

\address{H. Lange\\Mathematisches Institut,
              Universit\"at Erlangen-N\"urnberg\\Germany}
\email{lange@mi.uni-erlangen.de}

\address{R. E. Rodr{\'\i}guez\\Facultad de Matem\'aticas,
Pontificia Universidad Cat\'olica de Chile, Casilla 306-22,
Santiago, Chile}
\email{rubi@mat.puc.cl}

\address{A. M. Rojas \\Departamento de Matem\'aticas, Facultad de
Ciencias, Universidad de Chile, Santiago\\Chile}
\email{anirojas@uchile.cl }%

\thanks{The first, third and fourth author were supported by Fondecyt
grants 1060743, 1060742 and 11060468 respectively}%

\begin{abstract}
Let $X_1, \ldots, X_m$ denote smooth projective curves of genus
$g_i \geq 2$ over an algebraically closed field of characteristic
$0$ and let $n$ denote any integer at least equal to
$1+\max_{i=1}^m g_i$. We show that the product $JX_1 \times \cdots
\times JX_m$ of the corresponding Jacobian varieties admits the
structure of a Prym-Tyurin variety of exponent $n^{m-1}$. This
exponent is considerably smaller than the exponent of the
structure of a Prym-Tyurin variety known to exist for an arbitrary
principally polarized abelian variety. Moreover it is given by
explicit correspondences.
\end{abstract}
\maketitle


\section{Introduction}

A \textit{Prym-Tyurin variety of exponent} $q$ in a Jacobian $J$
is by definition an abelian subvariety of $J$ such that the
canonical polarization of $J$ induces the $q$-fold of a principal
polarization. According to a theorem of Matsusaka-Ran and 
Welters' criterion (see \cite[Section 12.2]{bl}), Prym-Tyurin varieties
of exponent $1$ are Jacobians. Welters showed in \cite{w} that,
roughly speaking, every Prym-Tyurin variety of exponent $2$ is a
classical Prym variety, that is the connected component
containing $0$ of the norm map $J\tilde{X} \ra JX$ of an \'etale
double covering $\tilde{X} \ra X$.

In \cite{M} Mumford showed that the product of Jacobians of two
general hyperelliptic curves occurs as a classical Prym variety.
In \cite{clrr} we generalized his construction by building
Prym-Tyurin varieties which are products of two Prym-Tyurin
varieties of smaller exponent. Here we generalize Mumford's result
in a different direction: We show that the product of an arbitrary
number of Jacobians of curves of genus $g$ occurs as a Prym-Tyurin
variety. To be more precise, the following theorem is a special
case of our main results.

\begin{thm} \label{thm1.1}
Let $X_1, \ldots, X_m$ be smooth projective curves of genus $g_i \geq 2$
over an algebraically closed field of characteristic $0$ and
$n \geq 1+ \max_{i=1}^m g_i$ an integer. Then the product
$$
JX_1 \times \cdots \times JX_m
$$
occurs as a Prym-Tyurin variety of exponent $n^{m-1}$ in a
Jacobian $J$ of dimension
$$
\dim J = n^{m-1}(\sum_{i=1}^m g_i + (m-1)n - m) + 1.
$$

\end{thm}

Our main results are Theorems \ref{prod} and \ref{thm4.4} below, which
also imply Mumford's above mentioned theorem as a special case.

For the proof we use a general result of \cite{clrr} which allows
the construction of new Prym-Tyurin varieties out of given ones.
We may assume that the base field is the field of complex numbers.
It is well known that any curve $X$ of genus $g$ admits simply
ramified coverings $X \ra \P^1$ of degree $n$, for all $n\geq
g+1$. We use any such covering to construct a structure of a
Prym-Tyurin variety of exponent $1$ on the Jacobian $JX$. Then one
only has to verify the hypotheses of the above mentioned result.

In Section 2 we recall the theory of \cite{clrr} for the special case
needed here. In Section 3 we work out the presentation of the Jacobians
as a Prym-Tyurin variety of exponent $1$. Section 4 contains the proofs
of our main results.

Throughout this paper, for any finite group $G$, $\chi_0$ will denote
its trivial representation. Also, for any subgroup $H$ of $G$,
$\rho_H^G$ will denote the representation of $G$ induced by the
trivial representation of $H$.

\section{Presentation of a Prym-Tyurin variety}

We want to apply the main result of \cite{clrr} to the symmetric group
$\mathbf{S}_n$ and its self-products. Since all rational irreducible
representations of these groups are absolutely irreducible, it suffices
to recall a special case of this theorem.\\

Let $G$ be a finite group, which later we assume to be $\mathbf{S}_n$.
Let $V_1, \ldots, V_r$ denote nontrivial, pairwise non-isomorphic,
absolutely irreducible rational representations of the same dimension
of the group $G$, and let $H$ be a subgroup of $G$ such that
for all $k = 1, \ldots, r$,
\begin{equation} \label{eqn2.1}
\dim V_k^H = 1 \; \text{ and } \; H \; \text{ is maximal
with this property}.
\end{equation}
Here \lq\lq maximal\rq\rq \, means that for every subgroup $N$ of $G$
with $H \subsetneqq N$ there  is an index $k$ such that $\dim V_k^N = 0$.

Choose a set of representatives
$$
\{ g_{ij} \in G \;| \; i = 1, \ldots, d \,
\text{ and } \, j = 1, \ldots, n_i \}
$$
for both the left cosets and right cosets of $H$ in $G$, and such that
$$
G = \sqcup_{i=1}^d Hg_{i1}H \, ,  \text{  and} \quad Hg_{i1}H =
\sqcup_{j=1}^{n_i} g_{ij}H = \sqcup_{j=1}^{n_i} Hg_{ij}
$$
are the decompositions of $G$ into double cosets, and of the double
cosets into right and left cosets of $H$ in $G$. Moreover, we assume
$g_{11} = 1_G$.\\

Now let $Z$ be a (smooth projective) curve with $G$-action and quotients
$$
\Pi: Z \ra \P^1 = Z/G \quad \text{and} \quad \pi : Z \to X := Z/H.
$$
In \cite{clrr} we defined a correspondence on $X$, which is given by
\begin{equation} \label{eqnD}
{\cD}(x) = \sum_{i=1}^d b_{i} \sum_{j=1}^{n_i} \pi g_{ij}(z).\
\end{equation}
for all $x \in X$ and $z \in Z$ with $\pi(z) = x$, where
\begin{equation} \label{eq:bs}
b_i := \sum_{k=1}^r \sum_{h \in H} \chi_{V_k}(hg_{i1}^{-1})
\end{equation}
is an integer for $i = 1, \ldots, d$. Moreover, denote
\begin{equation} \label{eqn2.3}
b := \gcd \{ b_1 - b_i \; | \; 2 \leq i \leq d \}.
\end{equation}
Let $\delta_{\cD}$ denote the endomorphism of the Jacobian
$JX$ associated to the correspondence ${\cD}$. We denote
by
$$
P_{\cD} := \im(\delta_{\cD})
$$
the image of the endomorphism $\delta_{\cD}$ in the
Jacobian $JX$ and call it the \textit{(generalized) Prym variety}
associated to the correspondence ${\cD}$.

Finally, let us recall the \textit{geometric signature} of a
Galois covering of curves with Galois group $G$. Let $C_1, \ldots,
C_t$ be pairwise different nontrivial conjugacy classes of cyclic
subgroups of $G$. Then this is by definition the tuple
$[\gamma,(C_1,s_1),\ldots, (C_t,s_t)]$, where $\gamma$ is the
genus of the quotient curve $Y$, the covering has a total of
$\sum_{j=1}^t s_j$ branch points in $Y$ and exactly $s_j$ of them
are of type $C_j$ for $j = 1, \ldots, t$; that is, the
corresponding points in its fiber have stabilizer belonging to
$C_j$.

Then \cite[Theorem 4.8]{clrr} can be stated as follows:

\begin{thm} \label{thm4.8}
Let $V_1, \ldots, V_r$ denote nontrivial pairwise non-isomorphic
absolutely irreducible rational representations of the group $G$
satisfying \eqref{eqn2.1} for a subgroup $H$ of $G$.

Suppose that the action of the finite group $G$ on a curve $Z$ has
geometric signature $[0;(C_1,s_1), \ldots , (C_{t},s_{t})]$
and satisfies
\begin{equation}  \label{eqn2.2}
\sum_{j=1}^t s_j \left( q \sum_{k=1}^r  (\dim V_k - \dim
V_k^{G_j}) - \left( [G:H] -|H \backslash G/G_j|\right) \right) = 0
\, ,
\end{equation}
where $G_j$ is a subgroup of $G$ of class $C_j$ and
\begin{equation} \label{eqq}
q = \frac{|G|}{b \cdot \dim V_1}.
\end{equation}
Then $P_{\cD}$ is a Prym-Tyurin variety of exponent $q$ in $JX$,
where $X = Z/H$.
\end{thm}

Furthermore, we showed in \cite[Section 4.4]{clrr} that
\begin{equation}
  \label{eq:dim}
\dim P_{\cD} = \sum_{i=1}^r [ - \dim V_1 +\frac{1}{2} \sum_{j=1}^t
s_j (\dim V_i - \dim V_i^{G_j}) ]
\end{equation}
and
\begin{equation}
  \label{eq:genus}
g_{X} = 1 -[G:H] + \frac{1}{2} \sum_{j=1}^t s_j ([G:H]
- |H\backslash G / G_j|).
\end{equation}

Observe that $b$, and hence $q$, depends only on the group $G$,
and not on its action on a given curve.

In the sequel we will use the following definition:  We say that
the construction of the Prym-Tyurin variety $P = P_{{\cD}}$ of
Theorem \ref{thm4.8} is a {\it presentation of} $P$ with respect
to the action of the group $G$, the subgroup $H$ and the set of
representations $\{V_1, \ldots, V_r\}.$

\section{The Jacobian of a curve as a Prym-Tyurin variety}

In this section we show that the Jacobian $JX$ of any curve $X$ of genus
$g \geq 2$ has a presentation as a Prym-Tyurin variety of exponent $1$
with respect to a symmetric group $\mathbf{S}_n$, a subgroup
$\mathbf{S}_{n-1}$ and the standard representation of $\mathbf{S}_n$.
First we need the following lemma.

\begin{lem}\label{pryms}
Let $G$ be a finite group, $V$ a nontrivial absolutely irreducible
rational representation of $G$, and $H$ a subgroup of $G$ such
that the representation $\rho_H^G$ of $G$ induced by the trivial
representation of $H$ has the isotypical decomposition
$$
\rho_H^G= \chi_0 \oplus V \, .
$$
Then the numbers $q$ and $b$, defined in equations \eqref{eqn2.3}
and \eqref{eqq} respectively, are given by
$$b=\frac{|G|}{\dim V}\quad \text{ and }\quad q=1.
$$
\end{lem}
\begin{proof}
The number of double cosets of $H$ in $G$ equals the following
character product
$$
|H\bs G/H| = \langle \rho_H^G,\rho_H^G\rangle_G \, ,
$$
and this is equal to two by our assumptions. Therefore we have two
double coset representatives $g_{11}=1_G$ and $g_{21}$. The number
$n_2$ of right cosets in $Hg_{21}H$ is $[G:H]-1=\dim V$. Using
$\langle \chi_0,V\rangle_G=0$, this gives
$$
0=\sum_{g\in G}\chi_{V}(g)=b_1+n_2b_2 \, ,
$$
where the $b_i$'s are given in \eqref{eq:bs}.

Moreover, using $b_1=|H|$ we obtain
$$
b_2=\frac{-b_1}{n_2}=\frac{-|H|}{\dim V}.
$$
By definition
$$
b=b_1-b_2=|H|+\frac{|H|}{\dim V}=
\frac{|H|(\dim V +1)}{\dim V}=\frac{|H|[G:H]}{\dim V} \, ,
$$
which implies both assertions.
\end{proof}

Now let $X$ be a smooth projective curve of genus $g\geq 2$ over an
algebraically closed field of characteristic 0. Recall that a covering
$f: X \ra \P^1$ of degree $n \geq 2$ is called \textit{simple} if the
fibre $f^{-1}(p)$  over every branch point $p \in \P^1$ consists of
exactly $n-1$ different points. According to \cite[Proposition 8.1]{fu},
$X$ admits a simple covering $f: X \ra \P^1$ of degree $n$ for any
$n \ge g+1$.

\begin{lem}  \label{lem3.2}
Let $f:X \ra \P^1$ denote a simple covering of degree $n \geq 2$.
Then

\medskip

\noindent {\em (a):} The covering $f$ does not factorize.

\noindent {\em (b):} The Galois group of the Galois closure $\Pi:
Z \ra \P^1$ of $f$ is the symmetric group $\mathbf{S}_n$.
\end{lem}

\begin{proof}
(a): Suppose $f = q \circ t$, with $t : X \ra X'$ and $q : X' \ra \P^1$.
If $q$ is nontrivial, then it has to be ramified as a covering of $\P^1$,
in which case $\deg t = 1$ since otherwise $f$ would not be a simple
covering.

(b): We may assume that the base field is the field of complex numbers so
that the monodromy group is well defined. Since the Galois group of the
covering $\Pi$ coincides with the monodromy group of the covering $f$, it
suffices to show that this is $\mathbf{S}_n$.

Since $X$ is irreducible and $f$ is simple, the monodromy group of $f$
is a transitive subgroup of $\mathbf{S}_n$ generated by transpositions.
By part (a), it is also primitive (or see \cite[Lemma 4.4.4]{serre}). It
is a well-known group theoretical theorem (see e.g. \cite[Satz 4.5,
p.171]{hu}) that any such subgroup coincides with the full group
$\mathbf{S}_n$.
\end{proof}

We consider the group
$$
G := \mathbf{S}_n
$$
as the permutation group of the set $\{1, \ldots, n\}$. As such,
it is generated by $\tau:=(1\;2)$ and
$\sigma:=(1\;2\;3\;\dots\;n)$. We denote by $Z$ the Galois closure
of (the simple covering) $f:X\to \P^1$. Therefore, $G$ acts on $Z$
with geometric signature $[0;(C_{\tau},s)]$, where $C_{\tau}$ is
the conjugacy class of any subgroup generated by a transposition,
and $s$ is the number of branch points of $f$.

Consider $H=\langle (1\;2\;3\;\dots\;n-1),(1\; 2) \rangle \simeq
\mathbf{S}_{n-1}$, the subgroup of $G$ fixing $n$, and let $V$ denote
the standard representation of $\mathbf{S}_n$ (of degree $n-1$) defined
by
$$
V: \; \sigma \rightarrow \left[ \begin{array}{ccccc}
0&& \cdots&0&-1\\
1&0&\cdots  &0&-1\\
0&1&\ldots & 0&-1\\
0 & 0 & \ddots\\
0&0 & \ldots  &1&-1\end{array} \right] \quad ,\quad \tau \rightarrow
\left[
\begin{array}{ccccc}
0&1&0 & \cdots&0 \\
1&0&0   &\cdots&0\\
0&0&1 &\cdots&0\\
0 & 0 & 0& \ddots & 0\\
0&0&\cdots & & 1\end{array} \right].
$$

In the following proposition we see that the Jacobian $JX$ of $X$ has a
presentation as a Prym-Tyurin variety of exponent $1$.

\begin{prop}
\label{prop3.3}
Let $X$ be a smooth projective curve of genus $g \geq 2$, and let $n$ be
any positive integer such that there exists an $n$-fold simple covering
$f : X \ra \mathbb{P}^1$. Let $G,H,V$, the curve $Z$ and the action of
$G$ on $Z$ be as above.

Then $X \simeq Z/H$ and the Jacobian $JX$ has a presentation as a
Prym-Tyurin variety of exponent $1$ with respect to the action of the
group $G$ on $Z$ with geometric signature $[0;(C_{\tau},s)]$, the subgroup
$H$ and the standard representation $V$ of $G$.
\end{prop}
\begin{proof}
First note that up to conjugacy $H$ is the only subgroup of index $n$
of $G$, therefore $X \simeq Z/H$. To prove that $JX$ has the desired
presentation, we use Theorem \ref{thm4.8}.

Certainly $\dim V^H = 1$ and $H$ is maximal with this property.
Moreover, it is easy to see that
$$
\rho_H^G=\chi_0\oplus V.
$$
Hence Lemma \ref{pryms} implies $q=1$. It remains to show that
equation \eqref{eqn2.2} is satisfied. Now we have
\begin{eqnarray*}
|H\bs G/\langle \tau \rangle|&=\langle \rho_H^G,
\rho_{\langle \tau \rangle}^G \rangle_G
=\frac{1}{2}\left[ (1+\chi_{V}(1_G))+(1+\chi_{V}(\tau))\right]=n-1
\end{eqnarray*}
and
$$
\dim V^{\langle \tau \rangle}=\langle V, \rho_{\langle \tau \rangle}^G
\rangle_G=\frac{1}{2}\left[ \chi_{V}(1_G)+\chi_{V}(\tau)\right]=n-2.
$$
where the middle equalities in both equations are due to Frobenius
reciprocity.
Inserting these values in equation \eqref{eqn2.2} completes the proof.
\end{proof}

\begin{rem}\label{rem:aes}

With the notation of Proposition \ref{prop3.3}, the correspondence
${\cD}$ on $X$ (defined in Equation \eqref{eqnD}) for this case is

$$
{\cD}(x) = \sum_{i=1}^2 a_{i} \sum_{j=1}^{n_i} \pi g_{ij}(z),
$$
where $n_1=1$, $n_2=n-1$, $a_1=(n-1)!$ and $a_2=-(n-2)!$. In fact,
the $g_{ij}$ may be chosen as $g_{11} = 1_G$, $g_{21}=(1 \; n)$,
$g_{22}= (2 \; n)$, $\ldots $ , $g_{n-1,n} = (n-1 \; n)$.

\end{rem}

\section{Products of Jacobians}

Fix integers $m \geq 2$ and $n \geq 2$. For each $i = 1, \ldots, m$, let
$f_i: X_i \ra \P^1$ denote a simple covering of degree $n$ with a smooth
projective curve $X_i$ of genus $g_i \geq 2$.

If $Z_i$ denotes the Galois closure of $f_i:X_i\to \P^1$, we have a
diagram for every $i = 1, \ldots, m$,
\begin{equation} \label{diag4.1}
\xymatrix{
      &  Z_i   \ar[dl]^{}_{\pi_i} \ar[dr]_{}^{\Pi_i}\\
     X_i  \ar[rr]^{f_i}^{}  & & \P^1.  \\
    }
\end{equation}
According to Lemma \ref{lem3.2} the Galois group of $\Pi_i$ is the group
$G := \mathbf{S}_n$. Moreover, the simple ramification of $f_i$ means
that $G$ acts with geometric signature $[0;(C_{\tau},s_i)]$ on $Z_i$,
where $C_{\tau}$ is the conjugacy class of the subgroups generated by a
transposition. The Riemann-Hurwitz formula implies
\begin{equation} \label{eqns_i}
s_i = 2(g_i +n - 1)
\end{equation} and thus
$$
g(Z_i) = \frac{n!}{2}(g_i +n -3) + 1.
$$
For every $i$, the curve $X_i \simeq Z_i/H$, with $H \simeq
\mathbf{S}_{n-1}$ of the previous section.

Consider the direct product group
$$
G^m := \times_{i=1}^m G \quad \text{ and its subgroup} \quad
H^m:= \times_{i=1}^m H \, ,
$$
and write $C_i$ for the conjugacy class in $G^m$ of the cyclic
subgroup generated by $\tau_i:=(1,\dots,\tau,1,\dots,1)$, with
$\tau = (1 \;2)$ in the $i$th coordinate.

\begin{lem}\label{lem:Z}
Suppose the branch loci of $f_i: X_i \ra \P^1$ are pairwise
disjoint in $\P^1$ (or, equivalently, the branch loci of $\Pi_i:
Z_i \ra \P^1$ are pairwise disjoint in $\P^1$). Then for all $m
\geq 1$ we have

\medskip

\noindent {\em (a):} The fibre product
$$
Z := Z_1 \times_{\P^1} Z_2 \times_{\P^1} \dots  \times_{\P^1} Z_m
$$
is a smooth projective curve of genus
$$
g(Z) =\frac{(n!)^m}{2}(\sum_{i=1}^mg_i + m(n-1)-2) + 1 ;
$$
it is a Galois covering  of $\P^1$ with geometric signature
$[0;(C_1,s_1), \ldots , (C_m, s_m)]$ and Galois group $G^m$.

\medskip

\noindent {\em (b):} The curve $X := Z/H^m$ coincides with the fibre product
$$
X = X_1 \times_{\P^1}  X_2 \times_{\P^1} \cdots \times_{\P^1} X_m
$$
and the genus of $X$ is

$$g(X) = n^{m-1}(\sum_{i=1}^m g_i + (m-1)n - m) + 1.$$
{\em (c):} The natural projections $q_i: X \ra X_i$ do not factorize
as

\begin{equation} \label{diag_ch}
\xymatrix{
        X   \ar[dr]^{}^{q_i^1} \ar[dd]_{}_{q_i}  &\\
     &\tilde{X_i}  \ar[dl]^{q_i^2}^{}   \\
 X_i.  \\
    }
\end{equation}
with $q_i^2$ nontrivial (cyclic) \'etale covering for $i = 1, \ldots, m$.
\end{lem}

\begin{proof}
Using induction one immediately checks that the fact that the
branch loci are disjoint implies that $Z$ is smooth. Classical
Galois theory implies that $Z$ is Galois over $\P^1$ with Galois
group $G^m$. The last statement of (a) is clear from the
definitions. It is a consequence of the universal property of the
fibre product that the curve $X:= Z/H^m$ is the fibre product over
$\P^1$ of all the $X_i$, $i=1,\dots m$. The genera are computed
using the Riemann-Hurwitz
formula. This completes the proof of (a) and (b).\\
(c): It suffices to consider $q_1$, so assume $q_1$ factorizes as
in diagram \eqref{diag_ch}. Since the covering $f_1:X_1\to \P^1$
is isomorphic to the covering $Z/H\times G^{m-1}\to \P^1$, it
follows that the curve $\tilde{X_1}$ is given by a quotient $Z/N$,
where $N$ is a subgroup of $G^m$ such that
\begin{equation} \label{eq:subn}
H^m \subseteq N \subseteq H \times G^{m-1}.
\end{equation}

The hypothesis that $q_1^2$ is \'etale implies that the covering
$\tilde{X_1}\to \P^1$ does not ramify over the branch points of
type $C_j$ for $j=2,\dots, m$. Therefore $N$ contains all these
conjugacy classes; since the class $C_j$ generates the
corresponding subgroup $\{1_G\}^{j-1} \times G \times
\{1_G\}^{m-j}$, we conclude that $H\times G^{m-1}\leq N$, which
together with \eqref{eq:subn} imply that $\tilde{X}_1 = X_1$.
\end{proof}

For $j = 1, \ldots, m$ we consider the following absolutely
irreducible rational representation of $G^m$
\begin{equation}\label{eqn:reps}
V_j:=\chi_0\otimes \dots \otimes \chi_0\otimes V \otimes
\chi_0\otimes \dots \otimes \chi_0,
\end{equation}
given by the outer tensor product of the trivial representation
$\chi_0$ ($m-1$ times) with the standard representation $V$ for $G$ in
the $j$-th component.

For $i=1, \ldots, m$, let $\cD_i \subset X_i \times X_i$ denote
the correspondence defined by \eqref{eqnD} with respect to
the action of the group $G$ on $Z_i$, the subgroup $H$ and
the representation $V$ of $G$.
Similarly let $\cD \subset Z \times Z$ denote the correspondence
defined by \eqref{eqnD} with respect to the action of the group
$G^m$ on $Z$, the subgroup $H^m$ and the representations
$V_1, \ldots, V_m$ of $G^m$. If $q_i: X \ra X_i$ denotes the
natural projection map, we have the following equality of divisors.

\begin{lem} \label{lem3.6}
$$
{\cD}=|H|^{m-1}\cdot \sum_{i=1}^m q_i^* {\cD}_i.
$$
\end{lem}

\begin{proof}
As in Section 2 we denote $d=|H\bs G/H|$ and
$\{g_{ij}:i=1,\ldots,d, \; j=1, \ldots ,n_i\}$. Therefore $|H^m\bs
G^m/H^m|=d^m$, $\{ g_{i_1 \ldots i_m} = (g_{i_1 1}, \dots ,g_{i_m
1}):i_k = 1,\ldots, d,\;  k=1,\dots,m\}$ are representatives for
the double cosets of $H^m$ in $G^m$, and $\{(g_{i_1j_1}, \dots
,g_{i_mj_m}):i_k = 1,\ldots, d,\; j_k=1, \ldots ,n_{i_k},
k=1,\dots,m\}$ are representatives of both left and right cosets
of the subgroup $H^m$ of $G^m$. According to \eqref{eqnD} as
applied to $G$, $H$ and $V$ we have
$$
{\cD}_{\nu}(x_{\nu}) = \sum_{i=1}^d a_{i} \sum_{j=1}^{n_i}
\pi_{\nu} g_{ij}(z_{\nu}) =: \sum_{i=1}^d a_{i} D_{\nu i}
$$
for all $x_{\nu} \in X_{\nu}$, where $z_{\nu} \in Z_{\nu}$ is a
preimage of $x_{\nu}$, $1\leq {\nu}\leq m$, and where
$$
a_i := \sum_{h \in H} \chi_{V}(hg_{i1}^{-1})
$$
are the same integers for all ${\cD}_{\nu}$.

Also note that according to \eqref{eqnD} applied to $G^m$, $H^m$
and $V_1 , \ldots , V_m$ we have
$$
\cD (x_1 , \ldots , x_m) = \sum_{1 \leq i_1 , \ldots , i_m \leq d}
b_{i_1 , \ldots , i_m} (D_{1i_1} , \ldots , D_{m i_m})
$$
where
\begin{align}
\notag b_{i_1 , \ldots , i_m} & = \sum_{k=1}^m \sum_{h \in H^m}
\chi_{V_k}(h g_{i_1 \ldots i_m}^{-1}) \\
& = \sum_{h \in H^m} \chi_{V_1}(h g_{i_1 \ldots i_m}^{-1})+ \ldots
+ \sum_{h \in H^m} \chi_{V_m}(h g_{i_1 \ldots i_m}^{-1}) \label{eq:ba} \\
\notag & = \sum_{h= (h_1 , \ldots , h_m) \in H^m} \chi_{V}(h_1
g_{i_1 1}^{-1}) + \ldots + \sum_{h= (h_1 , \ldots , h_m) \in H^m}
\chi_{V}(h_m g_{i_m 1}^{-1}) \\
\notag & = |H|^{m-1} a_{i_1} + \ldots + |H|^{m-1}a_{i_m} .
\end{align}

Therefore
\begin{align}
\notag \cD (x_1 , \ldots , x_m) & = |H|^{m-1} \sum_{1 \leq i_1 ,
\ldots , i_m \leq d} (a_{i_1} + \ldots + a_{i_m}) (D_{1i_1} ,
\ldots , D_{m
i_m}) \\
& = |H|^{m-1} ( \sum_{1 \leq i_1  \leq d} a_{i_1} \sum_{1 \leq i_2
, \ldots , i_m \leq d} (D_{1i_1} , \ldots , D_{m i_m}) + \ldots \label{eq:d}\\
\notag \mbox{ } & \mbox{ } \ \ \ \ldots + \sum_{1 \leq i_m  \leq
d} a_{i_m} \sum_{1 \leq i_1 , \ldots , i_{m-1} \leq d} (D_{1i_1} ,
\ldots , D_{m i_m}))
\end{align}

Now by definition we have $(q_{\nu}^*{\cD}_{\nu})(x_1,\dots,
x_m)=q_{\nu}^{-1}{ \cD}_{\nu}q_{\nu}(x_1,\dots, x_m)
=q_{\nu}^{-1}{\cD}_{\nu}(x_{\nu})$ for ${\nu}=1,\dots, m$.
Therefore
$$
(q_{\nu}^*{\cD}_{\nu})(x_1,\dots,
x_m)=q_{\nu}^{-1}(\sum_{i=1}^da_{i} D_{\nu i})= \sum_{i=1}^d a_{i}
q_{\nu}^{-1}(D_{\nu i})
$$
and we see from  \eqref{eq:d} that
$$
\cD (x_1 , \ldots , x_m) = |H|^{m-1}( q_1^{-1}(D_1(x_1)) + \ldots
+ q_m^{-1}(D_m(x_m))),
$$
from where the result follows.
\end{proof}

The following two theorems are the main result of the paper.

\begin{thm} \label{prod}
For each $i=1,\dots, m$, consider a simply ramified covering
$f_i:X_i\to \P^1$  of degree $n$, with $X_i$ of genus $g_i \geq
2$, and with pairwise disjoint branch loci. Let $\Pi_i:Z_i\to
\P^1$ be the Galois closure of $f_i$ over $\P^1$. Denote by $Z$
the fiber product of all the curves $Z_i$ over $\P^1$, and by $X$
the fiber product of the curves $X_i$ over $\P^1$.

Then the action of the group $(\mathbf{S}_n)^m$ on the curve $Z$
defines a Prym-Tyurin variety $P$ in the Jacobian $JX$, with $P$
of exponent ${q} = n^{m-1}$ and dimension
\begin{equation} \label{eqn4.7}
\dim P = \sum_{i=1}^{m}
g_i = \sum_{i=1}^m \dim JX_i.
\end{equation}
\end{thm}

\begin{proof}
As above we write $G = \mathbf{S}_n$ and consider the subgroup $H
= \langle (1\;2\;3\;\dots\;n-1),(1\; 2) \rangle \simeq S_{n-1}$ as
the stabilizer of a point in the general fibre of the map $f_i:
X_i \ra \P^1$ for $i = 1, \ldots,m$. Let again $V$ denote the
standard representation of $G$ and $V_i$ the representations of
$G^m$ defined in $\eqref{eqn:reps}$. First we observe that
equation \eqref{eqn2.1} is satisfied for the subgroup
$H^m:=(\mathbf{S}_{n-1})^{m}$ of $G^m$ and the representations
$V_1, \ldots, V_m$.

To see this, note that
$$
\dim (V_i)^{H^m}  = \langle V_i , \rho_{H^m}^{G^m} \rangle_{G^m}  =
\langle V , \rho_H^G \rangle_G = 1
$$
for all $i$. The maximality of $H^m$ with respect to this property
is a consequence of the fact that every $V_i$ occurs in
$\rho_{H^m}^{G^m}$.

Then we have to compute the exponent $q$ as defined by equation
\eqref{eqq} in this case; for this we need to compute the number
$b$ of equation \eqref{eqn2.3}. Using \eqref{eq:ba} of Lemma
\ref{lem3.6} we obtain
$$
{b_1}=m|H|^{m-1}a_1.
$$
The rest of the coefficients of $\cD$ are of the following form
$$
|H|^{m-1}((m-\ell)a_1+\ell a_2)) \;\text{ for }\; \ell=1,\dots m
$$
with $a_1$ and $a_2$ the coefficients of the correspondence
$\cD_i$ as in Remark \ref{rem:aes}. Therefore the differences are
of the following type
$$
|H|^{m-1}\ell(a_1-a_2) \;\;\text{ for }\; \ell=1,\dots m,
$$
where $a_1 - a_2$ are the corresponding numbers for the
correspondence $\cD_i$. According to Lemma \ref{pryms} we have
$a_1 - a_2 = \frac{|G|}{\dim V} = (n-2)! n, \; $ which implies
$$b = |H|^{m-1}(n-2)!n
$$ and hence
$$
q=\frac{|G^m|}{{b}\cdot \dim V} =
\frac{|\mathbf{S}_n|^m}{|\mathbf{S}_{n-1}|^{m-1}(n-2)!n
\cdot (n-1)} = n^{m-1}.
$$

Therefore the assertion follows from Theorem \ref{thm4.8} as soon
as we show that
\begin{equation} \label{eqn3.2}
\sum_{i=1}^m s_i \left( {q}\sum_{k=1}^m
(\dim V_k - \dim (V_k)^{\langle \tau_i\rangle})  -
([G^m:H^m] -|H^m \backslash G^m/\langle \tau_i\rangle|)\right) = 0.
\end{equation}
with $s_i$ given by \eqref{eqns_i} and where $\tau_i$ was defined just
before Lemma \ref{lem:Z}.

To see this, observe that
$$
  [G^m:H^m] = n^m \quad   \text{ and  } \quad
|H^m \backslash G^m/\langle \tau_i\rangle| = [G:H]^{m-1}
|H\backslash G /\langle \tau \rangle | = n^{m-1}(n-1).
$$
Moreover, we have
$$
\dim (V_k)^{\langle \tau_i \rangle}  = \langle V_k ,
\rho_{\langle \tau_i \rangle}^{G^m}\rangle_{G^m}  =
\left\{
  \begin{array}{ll}
     \langle V, \rho_{\langle \tau \rangle}^G \rangle_G = n-2,
     & \hbox{if $k = i$;} \\
     \dim   V = n-1 & \hbox{otherwise.}
  \end{array}
\right.
$$
Hence the left hand side of \eqref{eqn3.2} is equal to
$$
\sum_{i=1}^ms_in^{m-1}\left[(\dim V_i- \dim V_i^{\langle \tau_i \rangle})
 -1\right]=0.$$
Finally, the computation of the dimension is a consequence of
equation \eqref{eq:dim} using \eqref{eqns_i}.
\end{proof}

The fact that the Prym-Tyurin variety $P$ is constructed via a product of
groups suggests that it is a product itself. Moreover equation
\eqref{eqn4.7} indicates that it is the product of the Jacobian varieties
$JX_i$. The next theorem shows that this is in fact the case.

\begin{thm} \label{thm4.4}
Let the notation be as in Theorem \ref{prod}. Then the maps
$q_i: X \ra X_i$ induce an isomorphism
$$
JX_1 \times \dots  \times JX_m \stackrel{\simeq}{\lra} P
$$
of principally polarized abelian varieties.
\end{thm}
\begin{proof}
The map $q_1^* + \dots + q_m^*: JX_1 \times \dots \times JX_m \ra JX$
is an isogeny onto its image.
According to Lemma \ref{lem3.6} it maps $JX_1 \times \dots \times JX_m$
into $P$. From Theorem \ref{prod} we obtain that $q_1^* + \dots + q_m^*$
induces an isogeny
$$
JX_1 \times \dots \times JX_m \ra P.
$$
According to Lemmas \ref{lem3.2} and \ref{lem:Z}(c) the maps
$q_i:X \ra X_i$ do not factorize via a nontrivial
cyclic \'etale covering. From this we deduce, using
\cite[Proposition 11.4.3]{bl}, that the canonical polarization
of $JX$ induces a polarization of the same type on $P$ and
$JX_1 \times \dots \times JX_m$, namely the
$n^{m-1}$-fold of a principal polarization.
This implies that $q_1^* + \dots + q_m^*: JX_1 \times
\dots \times JX_m \ra P$ is an isomorphism.
\end{proof}

A first consequence of Theorems \ref{prod} and \ref{thm4.4}
is Theorem 1.1.

\begin{proof}[Proof of Theorem 1.1]
Let $X_1, \ldots, X_m$ be smooth projective curves of genus $g_i \geq 2$
for all $i$ and $n$ an integer at least equal to $1+\max_{i=1}^m g_i$.
According to \cite[Proposition 8.1]{fu} each $X_i$ admits a simple
covering $f_i: X_i \ra \P^1$ of degree $n$, since $n \geq g_i + 1$. If
necessary, we may move the branch points so that they become pairwise
disjoint. According to Lemma \ref{lem3.2} the Galois group of the Galois
closure $Z_i \ra \P^1$ is the symmetric group $\mathbf{S}_n$. Hence the
assumptions of Theorems \ref{prod} and \ref{thm4.4} are satisfied. The
formula for the dimension of $JX$ is a special case of \eqref{eq:genus}.
\end{proof}

The following corollary is a direct consequence of Theorems
\ref{prod} and \ref{thm4.4}. Note that Mumford's theorem
(see \cite[p.346]{M}) mentioned in the introduction is just the
special case $m=2$ of it, using Welters Theorem (see \cite{w})
which implies that these Prym-Tyurin varieties of exponent $2$ are
classical Prym varieties.

\begin{cor}
Let $X_1, \ldots, X_m$ denote hyperelliptic curves of genus
$g_i \geq 2$ for $i = 1, \ldots m$ whose hyperelliptic coverings
have pairwise disjoint ramification locus. Then the product
$$
JX_1 \times \cdots \times JX_m
$$
occurs as a Prym-Tyurin variety of exponent $2^{m-1}$ in a
Jacobian of dimension
$$
\dim JX = 2^{m-1}(\sum_{i=1}^m g_i + m-2) + 1.
$$
\end{cor}

\begin{rem}
An analogous corollary can be stated for any $k \geq 3$ in
the case of simply ramified  $k$-gonal curves
$X_1, \ldots,  X_m$. In particular, the case $m=2, k=3$ gives the
Prym-Tyurin varieties of exponent $3$ which were studied
in \cite{LRR} starting from a completely different geometric
set-up (see \cite[Theorems 4.1 and 4.2]{LRR}).
\end{rem}

\begin{rem}
It is well-known that any principally polarized abelian
variety of dimension $g$ occurs as a Prym-Tyurin variety of
exponent $2^{g-1}(g-1)!$ (see \cite[Corollary 12.2.4]{bl}
and use twice the principal polarization). Notice that
the exponent $n^{m-1}$ is considerably smaller than this number.
Moreover, here the Prym-Tyurin varieties are given by
an explicit correspondence whereas, in the general case,
they are given somewhat abstractly by successive hyperplane sections.
\end{rem}

\begin{rem}
Our method for constructing Prym-Tyurin structures on products of Jacobians also works for groups different from
$\mathbf{S}_n$. We only give one example, namely using the alternating group $\mathbf{A}_n$ 
and its standard representation. We omit the details of the proof of the following theorem, since they 
are completely analogous to the proof of Theorem \ref{thm1.1}.
\end{rem}
\begin{thm} \label{thm4.9}
Let $X_1 , \ldots, X_m$ denote general curves of 
genus $g_i \geq 3$ over an algebraically closed field of characteristic 0 and  and 
$n \geq 1 +2 \max_{i=1}^m g_i$ an integer.Then the
product
$$
JX_1 \times \cdots \times JX_m
$$
occurs as a Prym-Tyurin variety of exponent $n^{m-1}$ in a
Jacobian $J$ of dimension
$$
\dim J = 1+n^{m-1}(n(2m-1)-2m+2\sum_{i=1}^m g_i)
$$
as well as in a Jacobian of dimension
$$
\dim J = 1+n^{m-1}(n(m-1)-m+\sum_{i=1}^m g_i).
$$
\end{thm}

For the proof we use the following result, proven in \cite[Theorem 3.3]{MV}.
Let $g \geq 3$. Then a general curve of genus $g$ admits a cover to
$\P^1$ of degree $n$ with monodromy group $\mathbf{A}_n$  such that all
inertia groups are generated by a double transposition if and only if $n
\geq 2g + 1$.
The assertion also holds for three-cycles instead of double
transpositions.

So for each $n \geq 2g+1$ a
general curve $X$ of genus $g$ admits  two kind of coverings
$X \ra \P^1$. Both are of degree
$n$ and their corresponding Galois cover $\Pi : Z \ra
\P^1$ has as Galois group $\mathbf{A}_n$. In the first case the action of $\mathbf{A}_n$ on $Z$ has geometric
signature $[0, (C_1,s_1)]$, where $C_1$ is the conjugacy class in
$\mathbf{A}_n$ of the subgroup generated by $(1\, 2)(3\; 4)$, and in the second case geometric signature $[0,
(C_2,s_2)]$ where $C_2$ is the conjugacy class in $\mathbf{A}_n$ of the
subgroup generated by $(1\; 2\; 3)$.

As before, it may be proven that then the Jacobians $JX_i$ have a
presentation as a Prym-Tyurin variety with exponent $q=1$ with respect to
the group $\mathbf{A}_n$, the subgroup $\mathbf{A}_{n-1}$ and the standard
representation of $\mathbf{A}_n$. Then Theorem \ref{thm4.8} can be applied to complete the proof of Theorem \ref{thm4.9}.

Observe that the dimension of the Jacobian $J$ in the first case comes from the fact that
the geometric signature for the action on all the curves is $[0, (C_1,s_1)]$ and in the second case on all the curves is
$[0, (C_2,s_2)]$. Of course, one could also work out a mixed case.

\begin{rem}
One could also use non-absolutely irreducible representations for constructing Prym-Tyurin structures on products of Jacobians. 
The proof is essentially the same. Instead of Theorem \ref{thm4.8} one needs the more general theory of \cite{clrr}. 
\end{rem}

\bibliographystyle{amsplain}

\end{document}